\documentclass[12pt,reqno]{amsart}

\usepackage{multirow}
\usepackage{hhline}

\usepackage{mathtools}
\usepackage{times}
\usepackage[T1]{fontenc}
\usepackage{mathrsfs}
\usepackage{latexsym}
\usepackage[dvips]{graphics}
\usepackage[titletoc, title]{appendix}
\usepackage{epsfig}
\usepackage{amssymb}
\usepackage{amsmath,amsfonts,amsthm,amssymb,amscd}
\usepackage{color}
\usepackage{hyperref}
\usepackage{url}
\usepackage{breakurl}
\newcommand{\bburl}[1]{\textcolor{blue}{\url{#1}}}
\usepackage{comment}
\usepackage{listings}
\usepackage{float}
\usepackage{bbm, dsfont}

\newcommand{\be}{\begin{equation}}
\newcommand{\ee}{\end{equation}}
\newcommand{\seqnum}[1]{\href{https://urldefense.com/v3/__https://oeis.org/*1*7D*7B*5Crm__;IyUlJQ!!DZ3fjg!5yDnHeGSfGljEqbpuRw46hLMiV_jaL8rn2Yk7DqDl12SpLF1Wq4c5kckywJ55Db23mDqXkwmG3-qhLCpIIaQw7Y$   \underline{#1}}}

\newtheorem{thm}{Theorem}[section]

\newtheorem{cor}[thm]{Corollary}

\newtheorem{lem}[thm]{Lemma}

\newtheorem{exa}[thm]{Example}

\newtheorem{rek}[thm]{Remark}

\newcommand{\parentheses}[1]{\left( {#1}\right)}
\newcommand\numberthis{\addtocounter{equation}{1}\tag{\theequation}}

\numberwithin{equation}{section}

\newcommand{\ith}[1]{${#1}^{\text{th}}$}

\providecommand{\floor}[1]{\lfloor #1 \rfloor}

\newcommand{\expectation}[1]{\text{E}\left[#1\right]}
\newcommand{\probability}[1]{\text{Prob}\left( {#1}\right)}
\newcommand{\indicator}[1]{\mathbbm{1}_{\{{#1}\}}}

\begin{document}

\title[The German Tank Problem with Multiple Factories]{The German Tank Problem with Multiple Factories}

\author[Miller]{Steven J. Miller}
\email{\textcolor{blue}{\href{mailto:sjm1@williams.edu}{sjm1@williams.edu}},
\textcolor{blue}{\href{mailto:Steven.Miller.MC.96@aya.yale.edu}{Steven.Miller.MC.96@aya.yale.edu}}}
\address{Department of Mathematics and Statistics, Williams College, Williamstown, MA 01267, USA}

\author[Sharma]{Kishan Sharma}
\email{\textcolor{blue}{\href{mailto:kds43@cam.ac.uk}{kds43@cam.ac.uk}}, \textcolor{blue}{\href{mailto:sharmakd2002@gmail.com}{sharmakd2002@gmail.com}}}
\address{Emmanuel College, University of Cambridge, Cambridge, CB2 3AP, UK}

\author[Yang]{Andrew K. Yang}
\email{\textcolor{blue}{\href{mailto:aky30@cantab.ac.uk}{aky30@cantab.ac.uk}},
\textcolor{blue}{\href{mailto:andrewkelvinyang@gmail.com}{andrewkelvinyang@gmail.com}}}
\address{Emmanuel College, University of Cambridge, Cambridge, CB2 3AP, UK}

\begin{abstract}
During the Second World War, estimates of the number of tanks deployed by Germany were critically needed. The Allies adopted a successful statistical approach to estimate this information: assuming that the tanks are sequentially numbered starting from 1, if we observe $k$ tanks from an unknown total of $N$, then the best linear unbiased estimator for $N$ is $M(1+1/k)-1$ where $M$ is the maximum observed serial number. However, in many situations, the original German Tank Problem is insufficient, since typically there are $l>1$ factories, and tanks produced by different factories may have serial numbers in disjoint ranges that are often far separated.

In \cite{CGM}, Clark, Gonye and Miller presented an unbiased estimator for $N$ when the minimum serial number is unknown. Provided one identifies which samples correspond to which factory, one can then estimate each factory’s range and summing the sizes of these ranges yields an estimate for the rival’s total productivity. We construct an efficient procedure to estimate the total productivity and prove that it is effective when $\log l/\log k$ is sufficiently small. In the final section, we show that given information about the gaps, we can make an estimator that performs orders of magnitude better when we have a small number of samples.
\end{abstract}

\subjclass[2020]{62B05}

\keywords{German Tank Problem, sufficient statistics, complete statistics, sampling without replacement, linear estimators}

\thanks{We thank Cory Simons for conversations on considering multiple factories, which helped lead to this project. We also thank the referees at PUMP for their recommendations. With their help, we wrote this paper with the intention for it to be approachable to all mathematics students.}

\maketitle

\tableofcontents

\section{Introduction}

During World War II, the Germans had an advantage in their capacity to use tanks compared to the Allies. To appropriately respond to the Germans, the Allies tried estimating the number of tanks on the battlefields. To start, spies were used to attempt to gain information on the scale of German tank production. However, the values reported lay far off from the actual values.

During the war, the Allies noticed that captured tanks had serial numbers that they could use to their advantage as information. After the following assumptions were made:
\begin{enumerate}
    \item the serial numbers of tanks were consecutive positive integers,
    \item the first tank has a serial number of 1, and
    \item any tank is equally likely to be captured,
\end{enumerate}
people formulated an estimate for the total number of tanks produced, $N$, based on the serial numbers of the captured tanks. For a sample of $1\leq k\leq N$ captured tanks, and maximum observed serial number $M$, we have the unbiased estimate $\hat{N}$ for $N$ of
\begin{equation}
\label{eq:originalgtp}
    \hat{N} = M \parentheses{1+\frac{1}{k}} - 1.
\end{equation}
This proved to be a good estimate. See Lee and Miller \cite{LM} for a derivation. A comparison between the observed and estimated number of tanks in Figure \ref{fig:germantanknumbers} from Ruggles and Brodie \cite{Ru} shows that the estimate was much more effective than intelligence at estimating tank production. Later in Section \ref{sec:gtpum}, we show that $\hat{N}$ in (\ref{eq:originalgtp}) is in fact the unique minimum variance unbiased estimator.

\begin{figure}[h!]
\begin{equation*}
\begin{array}{c|c|c|c}
     \text{Month} & \text{Intelligence estimate} & \text{Statistical estimate} & \text{German records} \\
     \hline
     \text{June 1940} &  1,000 & 169 & 122 \\
     \text{June 1941} & 1,550 & 244 & 271 \\
     \text{August 1942} & 1,550 & 327  & 342
\end{array}
\end{equation*}
    \centering
    \caption{Statistics vs Intelligence estimates for German tank production.}
    \label{fig:germantanknumbers}
\end{figure}

This is a instance in which statistical inference can be applied to real world problems, and does a better job than a non-mathematical approach such as espionage. Kalu \cite{Ka} and Wikipedia \cite{Wi} go through a deeper history of the problem. See Prosdocimi \cite{P}, Simon \cite{S} and Cheng, Eck and Crawford \cite{CEC} for a sampling of other research in statistical estimation inspired by the German Tank Problem. 

This paper comes from investigating cases in which the first assumption is not valid. An example where it would not hold is if we suppose tanks were produced in \emph{multiple factories}. Suppose we have factories labelled $1,2,\dots,l$, and factory $i$ produces $N_i$ tanks that are labelled with consecutive positive integers. However, we now allow there to be a gap between the last serial number of factory $i$ and the first serial number of factory $i+1$. 

In the multiple factories problem, capturing tanks can be viewed mathematically as sampling uniformly \emph{without} replacement from the union of $l$ disjoint sets of consecutive positive integers with lengths $N_1,\dots,N_l$, or from the set

\begin{eqnarray}
T & = & \{1,\dots,N_1\} \ \cup \ \{N_1 \ +G_1 \ + \ 1,\dots,N_1 \ +G_1 \ +N_2\} \ \cup \nonumber\\& & \cdots \ \cup \ \left\{\sum_{i=1}^{l-1} (N_i \ + \ G_i) \ + \ 1,\dots, \sum_{i=1}^{l-1} (N_i \ + \ G_i) \ + \ N_l \right\},
\end{eqnarray}

\noindent where we define $G_i$ to be the \ith{i} gap, or the number of missing naturals between the \ith{i} and \ith{(i+1)} factories. We ask, given $x_1,\dots,x_k$ sampled uniformly without replacement from $T$, what inferences can we make about the total number of tanks produced, $N_{\text{tot}} = N_1 + \cdots +N_l$?

As in the original German Tank Problem (GTP), in the Multiple Factories Problem (MFP) we assume throughout that the smallest serial number is 1. Relaxing this assumption does not change much, but we chose not to for continuity with the old problem, and feel that it does not make the problem more interesting.

While this assumption tells us the minimum positive integer for the first factory, we do not know the minimum value of subsequent factories. This is a major difficulty, as the GTP estimator is only good if we know the minimum value. (If it is shifted away from 1 by a known constant, we can just shift back down our samples.) Fortunately, Clark, Gonye and Miller give an estimate for this exact problem in \cite{CGM}.

Clark et al. introduced the estimator

\begin{equation}
\label{eq:unknownMinGTP}
    \hat{N}_{\text{UM}} \ = S\left(1+\frac{2}{k-1}\right)-1
\end{equation}

\noindent for the size $N$ of a set of consecutive naturals $\{ a,a+1,\dots,a+N-1 \}$ when given $2\leq k\leq N$ samples $x_1,\dots,x_k$ that are sampled uniformly \emph{without} replacement. This estimator depends on the spread $S = M-W$, where $M = \max(x_1,\dots,x_k)$ and $W = \min(x_1,\dots,x_k)$. We henceforth refer to this statistical problem and estimator as the German Tank Problem with Unknown Minimum (GTP-UM). We show in Section \ref{sec:gtpum} that this is in fact the unique minimum variance unbiased estimator for this problem. 

We can now reduce the MFP to classifying which samples come from which factory, and then using the GTP to estimate the size of the first factory, and GTP-UM to estimate the size of subsequent factories.

\subsection{Main Results}

In Section \ref{sec:gtpum}, we extend the work of \cite{CGM} and derive further results about the German Tank Problem and German Tank Problem with Unknown Minimum.

We prove that the variance of the GTP-UM estimator is

\begin{equation}
    \text{\normalfont Var}(\hat{N}_{\text{UM}}) \ = \ \frac{2(N+1)(N-k)}{(k-1)(k+2)},
\end{equation}

\noindent and attain this closed form through rather remarkable identities on binomial coefficients. 

We also prove the new result that the GTP and the GTP-UM are the unique \emph{Minimum Variance Unbiased Estimators} (MVUEs) for their respective estimation problems.

\sloppy We then turn our attention to the Multiple Factories Problem. When our gaps $G_1,\dots,G_{l-1}$ are unknown, we argue in Section \ref{sec:motive} that to make any progress on the most general form of this problem, we must first be confident that we have a sample from every factory. We find that this happens with high probability when $\log l/\log k$ is sufficiently small, and prove a ``threshold effect'' for this happening with probability 1 in the asymptotic case.

In Section \ref{sec:firstapproach}, we find that provided we have enough samples, we can make a robust estimator for the total productivity, regardless of the factory lengths or gaps, giving simulations showing the test error of the estimator against the number of samples.

We run into difficulty when we do not have enough samples to be confident of having a sample from every factory. This is also likely to be the case in a practical situation. In Section \ref{sec:restriction}, we investigate the progress we can make after forcing the factory sizes to be equal and the gaps to be both equal and \emph{known}, that is $N_1 = \cdots = N_{l}$ and $G = G_1 = \cdots = G_{l-1}$. We find that even when we have a sample size smaller than the number of factories, we can still get estimates several orders of magnitude better than we do before assuming equal and known gap sizes.

\section{German Tank Problem with an Unknown Minimum}
\label{sec:gtpum}

\begin{thm} We have for the German Tank Problem with Unknown Minimum that

\begin{equation}\label{eqn:GTPUMvariance}
    \text{\normalfont Var}(\hat{N}_{\text{UM}}) \ = \ \frac{2(N+1)(N-k)}{(k-1)(k+2)}.
\end{equation}

\end{thm}

The key ingredients to these closed form results are the following rather remarkable binomial coefficient identities proved in Appendix \ref{sec:identity}. \\

\noindent Identity I: For all $N \geq k$:
\begin{equation*}\label{identityI}
    \sum_{m=k-b+1}^{N-b+1} m\frac{\binom{m-1}{k-b} \binom{N-m}{b-1}}{\binom{N}{k}} \ = \ \frac{\parentheses{N+1}\parentheses{k-b+1}}{\parentheses{k+1}}
\end{equation*}

\noindent and \\

\noindent Identity II: For all $N \geq k$:
\begin{equation*}\label{identityII}
\begin{split}
    \sum_{m=k-b+1}^{N-b+1} m^2\frac{\binom{m-1}{k-b} \binom{N-m}{b-1}}{\binom{N}{k}} \ & = \ \frac{(k-b+1)(k-b+2)(N+2)(N+1)}{(k+2)(k+1)} \\ &- \ \frac{\parentheses{N+1}\parentheses{k-b+1}}{\parentheses{k+1}}.
\end{split}
\end{equation*}

\begin{rek}
    We remark that \eqref{eqn:GTPUMvariance} compares to a variance of $\frac{(N+1)(N-k)}{k(k+2)}$ for the normal GTP, see \cite{LM}. This means that the variance of the GTP-UM is $\frac{2k}{k-1}$ times that of the GTP - or in other words, at least double. This is perhaps intuitive, since instead of uncertainties at just one side of the interval, we have uncertainties at both sides. 
\end{rek}

\begin{proof}

We want to find the variance of the spread $S$. From \cite{CGM}, we have that after observing $k$ uniform samples without replacement from $T = \{ a,a+1,\dots,a+N-1 \}$,

\begin{equation}
\text{Prob}(S=s) \ = \ \frac{(N-s)\binom{s-1}{k-2}}{\binom{N}{k}},
\end{equation}

\noindent for $s=k-1,k,\dots,N-1$, and zero otherwise. This is because there are $N-S$ choices for the minimum and maximum tanks, and $\binom{S-1}{k-2}$ ways to choose the other $k-2$ tanks in between. We can then calculate the variance of the spread by taking identities I and II with $b=2$:

\begin{equation}
\begin{split}
    \text{E}[S^2] \ & = \ \sum_{s=k-1}^{N-1} s^2 \cdot \frac{(N-s)\binom{s-1}{k-2}}{\binom{N}{k}} \\
    & = \ \frac{k(k-1)(N+1)(N+2)}{(k+1)(k+2)} \ - \ \frac{(N+1)(k-1)}{k+1}
\end{split}
\end{equation}

\noindent and

\begin{equation}
\begin{split}
    \left(\text{E}[S]\right)^2 \ & = \ \left(\sum_{s=k-1}^{N-1} s \cdot \frac{(N-s)\binom{s-1}{k-2}}{\binom{N}{k}}\right)^2 \\
    & = \ \left(\frac{(N+1)(k-1)}{k+1}\right)^2,
\end{split}
\end{equation}

\noindent so we have that

\begin{equation}
\begin{split}
    \text{Var}(S) \ &= \ \text{E}[S^2] \ - \ \left(\text{E}[S]\right)^2 \\
    &= \ \frac{(k-1)(N+1)}{(k+1)^2 (k+2)}\Big[(k+1)\big[k(N+2)-(k+2)\big]-(N+1)(k-1)(k+2)\Big] \\
    &= \ \frac{2(k-1)(N+1)(N-k)}{(k+1)^2 (k+2)}.
\end{split}    
\end{equation}

By the definition of $\hat{N}_{\text{UM}}$ we have that 

\begin{equation}
\label{eq:varspread}
    \text{Var}(\hat{N}_{\text{UM}}) \ = \ \frac{2(N+1)(N-k)}{(k+1)(k+2)}
\end{equation}

\noindent as required.

\end{proof}

\subsection{Proof GTP and GTP-UM are the Minimum-Variance Unbiased Estimators}
\label{sec:mvue}

The \emph{Lehmann-Scheff\'e theorem} states that any estimator that is unbiased for a given unknown quantity, and that depends on the samples only through a \emph{complete}, \emph{sufficient} statistic, is the \emph{unique} best unbiased estimator of that quantity - the minimum-variance unbiased estimator (MVUE). See \cite{CB} for details.

\begin{thm}

The GTP and the GTP-UM are the unique MVUEs for their respective estimation problems. 

\end{thm}

\begin{proof}

In \cite{CGM}, we see that these estimators are unbiased estimators of $N$. It remains to show that our estimators depend on the samples through only a complete, sufficient statistic.

For the GTP, our model is parameterised by $N$, and our estimator depends on the samples only through the maximum sample $M=\max(x_1,\dots,x_k)$. So if we can prove $M$ is a complete, sufficient statistic for $x_1,\dots,x_k$, we are done. Consider the distribution of the samples conditional on $M$. We have that the other $k-1$ samples are uniformly sampled without replacement from $\{ 1,\dots,M-1 \}$. Thus

\begin{equation}
\label{eqn:GTPpmf}
\begin{split}
    p(x_1,\dots,x_k | M, N) \ &= \ \binom{M-1}{k-1}^{-1} \\
    &= \ p(x_1,\dots,x_k | M),
\end{split}
\end{equation}

\noindent so conditional on $M$, our samples do not depend on $N$. Thus $M$ is a sufficient statistic by definition.

To show $M$ is complete, we need to show that for every measurable function $g$, if $\text{E}_N [g(M)] = 0$ for all $N\geq k\geq 1$, then $\text{Prob}_N (g(M) = 0) = 1$ for all $N\geq k\geq 1$. Indeed, using the probability mass function \eqref{eqn:GTPpmf}, we yield

\begin{equation}
\label{eqn:completeGTP}
    \text{E}_N [g(M)] = \binom{N}{k}^{-1} \ \sum_{m=k}^N g(m) \binom{m-1}{k-1}.
\end{equation}

Suppose $\text{E}_N [g(M)] = 0$ for all $N\geq k\geq 1$. By considering $N=1$ in (\ref{eqn:completeGTP}), we have that $g(1)=0$. By considering $N=2$, we then must have that $g(2)=0$, and so on. Since $M$ only takes values in $\{1,\dots,N\}$, we have that $g(M) = 0$ for all $N\geq k\geq 1$. Thus $M$ is a complete statistic.

For the GTP-UM, our model is parameterised by $(a,N)$, and our estimator depends on the samples only through the maximum sample $M=\max(x_1,\dots,x_k)$ and minimum sample $W= \min(x_1,\dots,x_k)$. So if we can prove $(W,M)$ is a complete, sufficient statistic for $x_1,\dots,x_k$, we are done. Consider the distribution of the samples conditional on $(W,M)$. We have that the other $k-2$ samples are uniformly sampled without replacement from $\{ W+1,\dots,M-1 \}$. Thus

\begin{equation}
\label{eqn:GTPUMpmf}
\begin{split}
    p(x_1,\dots,x_k | W,M,a,N) \ &= \ \binom{M-W-1}{k-2}^{-1} \\
    &= \ p(x_1,\dots,x_k | W,M),
\end{split}
\end{equation}

\noindent so conditional on $(W,M)$, our samples do not depend on $(a,N)$. Thus $(W,M)$ is a sufficient statistic by definition.

To show $(W,M)$ is complete, we need to show that for every measurable function $g$, if $\text{E}_{a,N} [g(W,M)] = 0$ for all positive integers $a$ and $N\geq k\geq 2$, then $\text{Prob}_{a,N} (g(W,M) = 0) = 1$ for all positive integers $a$ and $N\geq k\geq 2$. Indeed, using the probability mass function \eqref{eqn:GTPUMpmf}, we yield

\begin{equation}
\label{eqn:completeGTPUM}
    \text{E}_{a,N} [g(W,M)] \ = \ \binom{N}{k}^{-1} \ \ \sum_{w=a}^{a+N-k} \sum_{m=a+k-1}^{a+N-1} g(m) \binom{m-w-1}{k-2}.
\end{equation}

Suppose $\text{E}_{a,N} [g(W,M)] = 0$ for all positive integers $a$ and $N\geq k\geq 2$. By considering $N=2$ in (\ref{eqn:completeGTPUM}), we have that $g(a,a+1)=0$ for all positive integers $a$. By considering $N=3$, we then must have that $g(a,a+2)=0$ for all positive integers $a$, and so on. Since $W,M$ only take values in $\{a,\dots,a+N-1\}$ with $M>W$, we have that $g(W,M) = 0$ for all positive integers $a$ and $N\geq k\geq 2$. Thus $(W,M)$ is a complete statistic. \end{proof}

\section{A First Approach at the Multiple Factories Problem}
\label{sec:motive}

To motivate this section, consider an example with small numbers. 

\begin{exa}

Suppose we know that there are two factories of equal size, and obtain some small sample, say, $\{5,22,114,124\}$. For all we know, all four samples are from the first factory (at least a 6\% chance\footnote{Since $\frac{\binom{N}{4}}{\binom{2N}{4}}\geq 0.06$ for $N\geq 124$.}). Indeed, if we knew the gap between the factories to be large, say, 100, we are certain of every sample being from the first factory. In this situation we could then estimate its size to be around 150, and thus estimate the total number of tanks to be double that, around 300.

However, when the gap between the factories is unknown, it could well be, say, 75. Conditional on this, it is more likely that our two largest samples are from the second factory, and each of our factories have size around 30 (i.e., $S = \{1,\dots,30\}\cup \{105,\dots,135\}$), giving a total number of tanks of around 60! This is a huge amount of variance we cannot control.

\end{exa}

In general, we can define a reordering $x_{(1)},\dots,x_{(k)}$ of our sample $x_1,\dots ,x_k$ such that $x_{(1)}<\dots<x_{(k)}$. Suppose for now we still know nothing about the gaps $G_1,\dots,$ $G_{l-1}$, but we make the assumption that we do have a sample from all of the factories. In Section \ref{sec:probmiss} we find the probability this assumption is valid. 

Our assumption means that each gap between any consecutive factories is contained inside a gap between some two consecutive samples $x_{(i)}$ and $x_{(i+1)}$. Regardless of gap size, these gaps are more likely to be found in the larger gaps between consecutive samples. One natural idea is to consider the set of gaps $\{x_{(2)}-x_{(1)},\dots,x_{(k)}-x_{(k-1)}\}$ and take the $l-1$ largest members, and split our sample into $l$ sub-samples $X_1,\dots,X_l$ at these gaps. 

This risks choosing the wrong gaps. However, if the factories are spaced far apart, and we have a sample from every factory, this is unlikely to happen. And if the factories are not spaced far apart, then the estimation will not be too far off even if the wrong gaps are chosen. We expand this idea into a full approach at the MFP in Section \ref{sec:firstapproach} and evaluate its performance.

\subsection{Probability of missing a factory}
\label{sec:probmiss}

We now find the probability $P_{N,l,k}$ that $k$ samples, chosen uniformly without replacement from $l$ factories, all of size $N$, are missing a sample from one or more factories. We make all factories to be the same size as if not, smaller factories are more likely to be missing, but are also less significant if unaccounted for.

Let us index our factories with the set $\{1,\dots,l\}$. Note that we must have $k\leq Nl$. Let $e_T$ be the event we are missing samples from factories in $T \subseteq \{1,\dots,l\}$. Note that this definition is for a general set $T$ of factories.

\begin{thm}\label{thm:PNLKformula}We have

\begin{equation*}
    P_{N,l,k} \ = \ {Nl \choose k}^{-1} \ \sum_{i=1}^{l-1} (-1)^{i+1} \ {l \choose i} \ {N(l-i) \choose k}.
\end{equation*}

\noindent where $P_{N,l,k}$ is defined as above.

\end{thm}

\begin{proof}

We have

\begin{equation}
    \text{Prob}(e_T) \ = \ \frac{{N(l-|T|) \choose k}}{{Nl \choose k}}
\end{equation}

\noindent by counting the number of ways we can select $k$ samples from $l-|T|$ factories. Observe there are ${l \choose i}$ choices of $T$ with $|T|=i$. By inclusion-exclusion, we have

\begin{equation}
\begin{split}
    P_{N,l,k} \ &= \ \sum_{i=1}^{l} (-1)^{i+1} \ {l \choose i} \ \text{Prob}(e_T,\text{ where } |T|=i) \\
    &= \ {Nl \choose k}^{-1} \ \sum_{i=1}^{l-1} (-1)^{i+1} \ {l \choose i} \ {N(l-i) \choose k}
\end{split}
\end{equation}

\noindent as required.

\end{proof}

\begin{figure}[h]
\begin{center}
\scalebox{.57}{\includegraphics{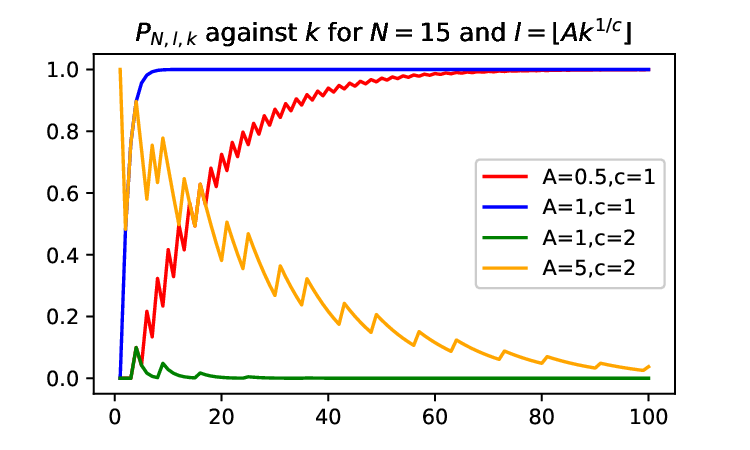}}\ \scalebox{.57}{\includegraphics{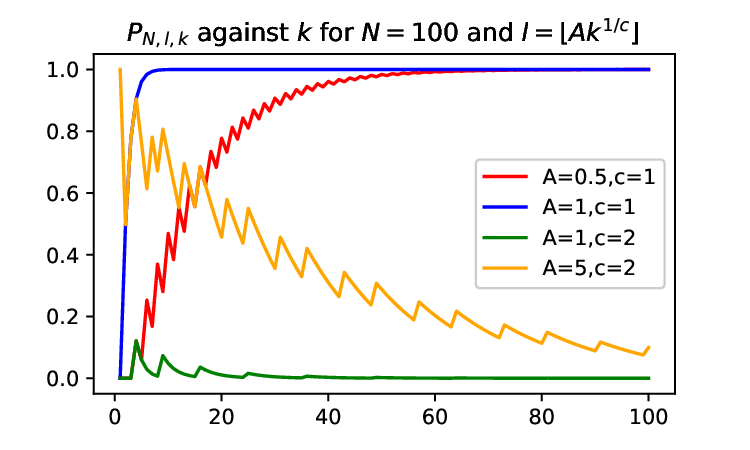}}
\caption{\label{figure:probmiss} Probability of missing at least one factory $P_{N,l,k}$ against the number of samples $k$ for different values of factory size $N$ and number of factories $l$. The zagged behaviour occurs due to the rounding of $l$ to integer values.}
\end{center}
\end{figure}

In Figure \ref{figure:probmiss} we see that the size of $l$ relative to $k$ is critical to whether the probability $P_{N,l,k}$ is closer to 1 or 0. When we set $l=\lfloor Ak^{1/c}\rfloor$, we observe that for $c=1$, we are increasingly likely to miss a factory as $k$ increases, but when $c=2$, we are increasingly likely to \emph{not} miss any factories. This is reasonable as when we have a lot of tanks relative to the number of factories we are very likely to have at least one from each.

We choose such an expression for $l$ because we require $l$ to be an integer, and in the next section, we show that when $N\to \infty$, the limiting behaviour of $P_{N,l,k}$ as $k,l\to \infty$ is dependent on the relative speed $l$ grows relative to $k$, or $\frac{\log l}{\log k}\to\frac{1}{c}$.

\subsection{Limiting Behaviour of $P$ as $N,l,k \to \infty$}

An obvious issue with the work in the previous section is that our probability $P_{N,l,k}$ depends on $N$, which is unknown. In this section, we investigate the limiting behaviour of $P_{N,l,k}$. Intuitively, taking $N \to \infty$ results in sampling without replacement converging to i.i.d. uniform sampling. Indeed,

\begin{lem}
\label{lem:probe_T}
As $N \to \infty$, for fixed $l,k$, 

\begin{equation*}
    \text{\emph{Prob}}(e_T,\text{ where } |T|=i) \ =  \ \frac{{N(l-i) \choose k}}{{Nl \choose k}} \ \to \ \left(\frac{l-i}{l}\right)^k.
\end{equation*}
\end{lem} 

\begin{proof} We have, as $N\to\infty$,

\begin{equation}
\begin{split}
    \frac{{N(l-i) \choose k}}{{Nl \choose k}} \ &= \ \frac{N(l-i)}{Nl} \cdot \frac{N(l-i)-1}{Nl-1} \cdots \frac{N(l-i)-k+1}{Nl-k+1} \\
    & \to \ \left(\frac{l-i}{l}\right)^k
\end{split}
\end{equation}

\noindent as required.

\end{proof}

Note that we can apply Lemma \ref{lem:probe_T} to Theorem \ref{thm:PNLKformula} to get

\begin{cor}
\label{cor:limN}
    As $N \to \infty$, for fixed $l,k$,

\begin{equation*}
    P_{N,l,k} \ \to \ \sum_{i=1}^{l-1} (-1)^{i+1} \ {l \choose i} \ \left(\frac{l-i}{l}\right)^k.
\end{equation*}
\end{cor}

This expression is nice as it is independent of $N$. However, inspired by the limiting behaviour exhibited in Figure \ref{figure:probmiss}, we can go further and analyse the limiting behaviour of $P_{N,l,k}$ when also taking $k,l \to \infty$ at different speeds. We strike upon an interesting result, exhibiting an example of a ``threshold effect''.

\begin{cor}
    As $N,k \to \infty$, for fixed $l$, independent of rates,

\begin{equation*}
    P_{N,l,k} \ \to \ 0.
\end{equation*}
\end{cor}

\begin{proof}
    Taking $k \to \infty$ in Corollary \ref{cor:limN} yields the desired result.
\end{proof}

\begin{rek}
  In other words, we have a sample from every factory \emph{almost surely} when $N,k \to \infty$.  
\end{rek}

We now investigate limiting behaviour as $N, l, k \to \infty$, more specifically how differing growth rates of 

\begin{thm}
\label{thm:almostsure}
    Let $l=\lfloor A k^c \rfloor$ where $A,c>0$ are constants. Taking $N,k \to \infty$, (or equivalently $N,l\to\infty$) we have:
\begin{itemize}
    \item if $c\geq1$, then $P_{N,l,k}\to 1$ and we are missing a factory almost surely, and
    \item if $c<1$, then $P_{N,l,k}\to 0$ and we have a sample from every factory almost surely.
\end{itemize}
\end{thm}

\begin{proof}
    Define $\indicator{i}$ as below,
    
    \begin{equation}
        \indicator{i} \ = 
        \begin{cases}
            1 &\text{if factory $i$ has no sample in it} \\
            0 &\text{if factory $i$ has at least one sample in it.}            
        \end{cases}
    \end{equation}
    
    Then using \eqref{lem:probe_T} with $|T| \ = \ |\{i\}| \ = \ 1$, as $N \to \infty$, for fixed $l,k$, 
    
    \begin{equation}
        \begin{split}
            \expectation{\indicator{i}} \ &= \ \probability{\indicator{i} \ = \ 1}\\
            &= \ \left(1 - \frac{1}{l}\right)^k. \\
        \end{split}
    \end{equation}

    Further, using \eqref{lem:probe_T} with $|T| \ =  \ |\{i,j\}| \ = \ 2$, we find that for $i \neq j$,

    \begin{equation}
        \begin{split}
            \expectation{\indicator{i}\indicator{j}} \ &= \ \probability{\indicator{i}\indicator{j} \ = \ 1}\\
            &= \parentheses{1 - \frac{2}{l}}^k,
        \end{split}
    \end{equation}

    \noindent and for $i = j$, $\indicator{i}^2 \ = \ \indicator{i}$, so

    \begin{equation}
        \begin{split}
            \expectation{\indicator{i}^2} \ &= \ \parentheses{1 - \frac{1}{l}}^k.
        \end{split}
    \end{equation}

    By linearity of expectation, if $Q_{l,k}$ is the proportion of factories with no sample in them (as $N \to \infty$), we obtain 

    \begin{equation} \label{eq: Q expectation formula}
        \begin{split}
            \expectation{Q_{l,k}} \ &= \ \expectation{\frac{1}{l} \sum_{i=1}^l \indicator{i}} \\
            &= \ \frac{1}{l} \sum_{i=1}^l \expectation{\indicator{i}} \\
            &= \ \parentheses{1 - \frac{1}{l}}^k. \\
        \end{split}
    \end{equation}

    Similarly for the variance, we obtain

    \begin{equation}
        \begin{split} \label{eq: Q variance formula}
            \text{Var}\parentheses{Q_{l,k}} \ &= \ \expectation{Q_{l,k}^2} - \expectation{Q_{l,k}}^2 \\
            &= \ \frac{1}{l^2} \expectation{\parentheses{\sum_{i=1}^l \indicator{i}}^2 } - \parentheses{1 - \frac{1}{l}}^{2k}\\
            &= \ \frac{1}{l^2} \sum_{i=1}^l \expectation{\indicator{i}^2} + \frac{1}{l^2}\sum_{i \neq j} \expectation{\indicator{i} \indicator{j}} - \parentheses{1 - \frac{1}{l}}^{2k}\\
            &= \frac{1}{l} \parentheses{1 - \frac{1}{l}}^k + \parentheses{1 - \frac{1}{l}} \parentheses{1 - \frac{2}{l}}^k - \parentheses{1 - \frac{1}{l}}^{2k}\\
            &= \ \underbrace{\frac{1}{l} \parentheses{\parentheses{1 - \frac{1}{l}}^k -\parentheses{1 - \frac{2}{l}}^k}}_{:= \rho} + \underbrace{\parentheses{\parentheses{1 - \frac{2}{l}}^k -  \parentheses{1 - \frac{1}{l}}^{2k}}}_{:= \sigma}. 
        \end{split}
    \end{equation}

    Observe that $|\rho|$ is bounded from above by $2l^{-1}$ and by the squeeze theorem we get that $\rho \to 0$ as $l \to \infty$. Next consider  $k = \alpha l^{{1}/{c}}$ where $\alpha = A^{-1}$. We ignore the floor function as we can always bound the floor function $l=\lfloor A k^c \rfloor$ from above and below in the limit by some continuous functions of the form $\hat{l} = B k^{c}$. We want to show that $\sigma \to 0 \text{ as } l,k \to \infty \ \ \text{for all } A,c > 0$, so that we can apply Chebyshev's inequality \cite{Mil}  as the variance will tend to 0.

     We now look at the possible cases.    
    \begin{itemize}
        \item $c > 1$: Using \eqref{eq: Q expectation formula}, we obtain,

        \begin{equation}
            \expectation{Q_{l,k}} \ = \ \parentheses{1 - \frac{1}{l}}^{\alpha l^{1/c}} .
        \end{equation}

        In particular for all $\gamma > 0$, there exists $L > 0$ such that if $l  > L$ then $l^{1/c} \leq \gamma l$. Since $\parentheses{1 - 1/l} < 1$, we have that 

        \begin{equation}
            \parentheses{1-\frac{1}{l}}^{\alpha l^{1/c}} \geq \parentheses{1 - \frac{1}{l}}^{\alpha \gamma l} \to e^{-\alpha \gamma} \text{ as } l \to \infty.
        \end{equation}

        As $\gamma \to 0$, $e^{- \alpha \gamma} \to 1$, so since $Q_{l,k}$ is bounded from above, we have that $\expectation{Q_{l,k}} \ \to \ 1$ as $l \to \infty$. In this case, we do not need to calculate $\sigma$ since $Q_{l,k}$ is bounded from above by $1$ and its expectation tends to $1$, therefore it tends to 1 in probability.  


        \item $c = 1$: We have $k = \alpha l$. Substituting this into \eqref{eq: Q expectation formula} yields
        \begin{equation}
            \begin{split}
                \expectation{Q_{l,k}} \ &= \ \parentheses{1 - \frac{1}{l}}^{\alpha l} \\
                &\to \ e^{-\alpha} \text{ as } l \to \infty. 
            \end{split}
        \end{equation}

        We continue in looking at $\sigma$ in \eqref{eq: Q variance formula}. Observe that 

        \begin{equation}
            \parentheses{1 - \frac{2}{l}}^{\alpha l} \ \to \ 1 - e^{-2 \alpha} \text{ as } l \to \infty
        \end{equation}

        \noindent and 

        \begin{equation}
            \parentheses{1 - \frac{1}{l}}^{2\alpha l} \ \to \ 1 - e^{-2 \alpha} \text{ as } l \to \infty.
        \end{equation}

        Therefore, taking the difference of these, 

        \begin{equation}
            \sigma \ \to \ 0 \text{ as } l \ \to \  \infty.
        \end{equation}

        Applying Chebyshev's inequality, for all $ \epsilon > 0$,

        \begin{equation}
            \probability{|Q_{l,k} - e^{-\alpha}| < \epsilon} \ \to \ 0 \text{ as } l \to \infty.
        \end{equation}

        Consequently, we expect that we are going to miss a factory and in fact that in the limit we expect a proportion of $e^{-\alpha}$ factories to have no samples.

        \item $c < 1$: We have $k = \alpha l^{{1}/{c}}$. For all $\gamma > 0$, there exists $L > 0$ such that if $l  > L$ then $l^{{1}/{c}} \geq \gamma l$. Since $\parentheses{1 - 1/l} < 1$, we have that 

        \begin{equation}
            \parentheses{1-\frac{1}{l}}^{\alpha l^{1/c}} \ \leq \ \parentheses{1 - \frac{1}{l}}^{\alpha \gamma l} \to e^{-\alpha \gamma} \text{ as } l \to \infty.
        \end{equation}

        As $\gamma \to \infty$, we can see that 

        \begin{equation}
            \expectation{Q_{l,k}} = \parentheses{1-\frac{1}{l}}^{\alpha l^{1/c}} \to 0 \text{ as } l \to \infty.
        \end{equation}

        Observe that each term in $\sigma$ tends to $0$ by the same reasoning as above as $l \to \infty$. Therefore, applying Chebyshev's inequality again,

        \begin{equation}
            \probability{|Q_{l,k}| < \epsilon} \ \to \ 0 \text{ as } l \to \infty.
        \end{equation}
        
    \end{itemize}\end{proof}

\begin{figure}[h]
\begin{center}
\scalebox{.57}{\includegraphics{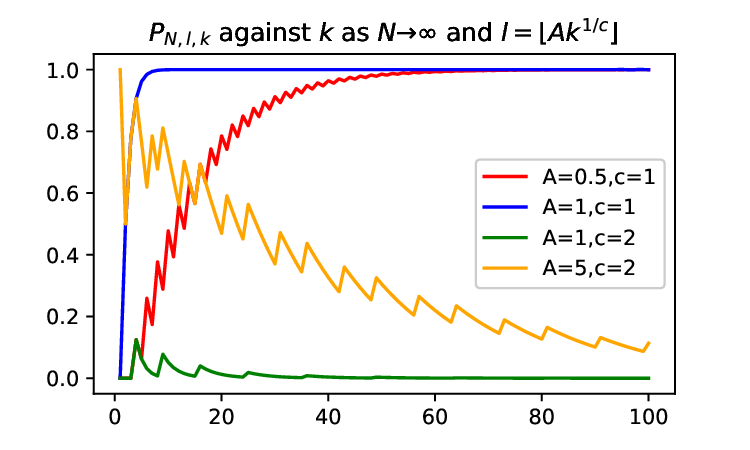}}\ \scalebox{.57}{\includegraphics{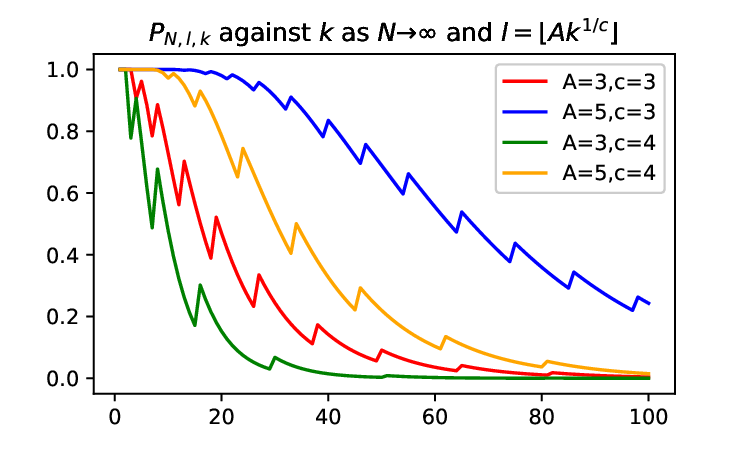}}
\caption{\label{figure:limNprobmiss} Probability of missing at least one factory $P_{N,l,k}$ against the number of samples $k$ for different values of the number of factories $l$ and $N\to\infty$.}
\end{center}
\end{figure}

We conclude that the rate of growth of $l$ relative to $k$ is the critical factor to whether we can be confident of having a sample from every factory, even asymptotically. Figure \ref{figure:limNprobmiss} gives empirical evidence of this behaviour that we have proved.

\subsection{A First Approach using GTP and GTP-UM}
\label{sec:firstapproach}

Recall that for the MFP we assume the number of factories is a known positive integer $l>1$, and each have sizes $N_1\dots,N_l$. We aim to estimate $N_{\text{tot}} = N_1 + \cdots + N_l$ from a uniform without replacement sample of tanks $x_1,\dots x_k$. We can order our sample such that $x_{(1)}<\dots<x_{(k)}$. We assume that each factory has contributed at least one sample, and the validity of this assumption is discussed in the previous sub-section.

At the start of the section, we motivated an approach to split our sample into $l$ sub-samples $X_1,\dots,X_l$ at the points of the $l-1$ largest gaps between consecutive samples. We can then apply the GTP-UM to estimate the size of each of the factories corresponding to our sub-samples $X_2,\dots,X_l$. We can do better with the first factory, and use the original GTP, since the minimum is known to be 1. We proved in Section \ref{sec:gtpum} that GTP-UM has variance at least double that of the GTP.

The idea is then to make an estimate $\hat{N}_i$ for $N_i$ using $X_i$, and estimate $N_{\text{tot}}$ by $\hat{N}_{\text{tot}} = \hat{N}_1 + \cdots + \hat{N}_l$.

There are two immediate questions with this approach that we will resolve somewhat crudely. One question is that GTP-UM requires at least two samples to make an estimation. When we have $|X_i|=1$ for some $i\neq1$, this is a ``bad'' sub-sample, and we must try and infer the size of the corresponding factory in other ways.

One reasonable approach is to estimate the size of the factory using the (estimated) sizes of the other factories and the number of samples. Let $I=\{i \ : \ i=1 \text{ or } |X_i|>1\}$, i.e., an index set for our ``good'' sub-samples. Suppose we have estimates $\{ \hat{N}_i \ : i\in I \}$. Let $\mathcal{N}=\sum_{i\in I} \hat{N}_i$, and let $\tilde{k}=\sum_{i\in I} |X_i|$.

We can now make a rough estimate for the size of the singleton sample factories ($|X_i|=1$ for some $i\neq1$): 

\begin{equation}
\label{eq:singlesampleestimator}
\hat{N}_{\text{bad}} \ = \ \frac{\mathcal{N}}{\tilde{k}}.
\end{equation}

We are using the fact that our samples are taken uniformly to say that the number of samples that came from a particular factory is roughly proportional to the size of the factory. The more singleton factories we have relative to non-singleton factories, the less confidence we have in this estimate due to the low sample size.

The second immediate question is resolving ties for the ${(l-1)}^{\text{st}}$ largest gap. Since we have little to no information to go on to decide which gaps to choose, we will just resolve ties by choosing the leftmost gaps in the tie. 

\begin{rek}
    Perhaps one improvement in the situation of ties could be choosing the largest gaps such that they prevent the formation of ``bad'' sub-samples if possible. However, note that when we have factory sizes that are reasonably large compared to the number of samples, ties are unlikely.
\end{rek}

MATLAB code that implements this method and simulates the MFP with independent samples and estimates the mean squared error of the estimator can be found in Appendix. \ref{sec:codeandy}.

\subsection{Simulated Performance of MFP Estimation}

\begin{figure}[h]
\begin{center}
\scalebox{.45}{\includegraphics{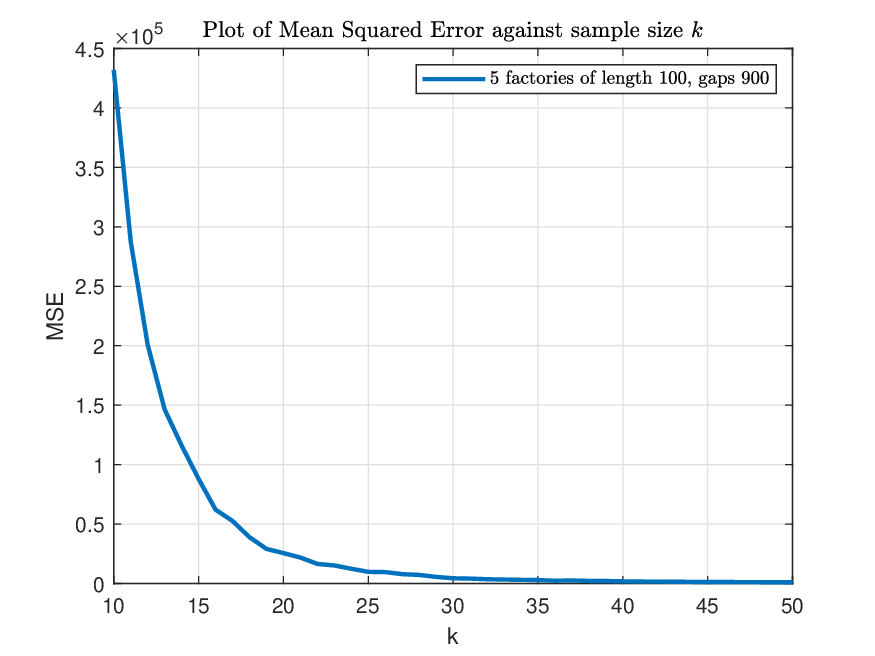}}\ \scalebox{.45}{\includegraphics{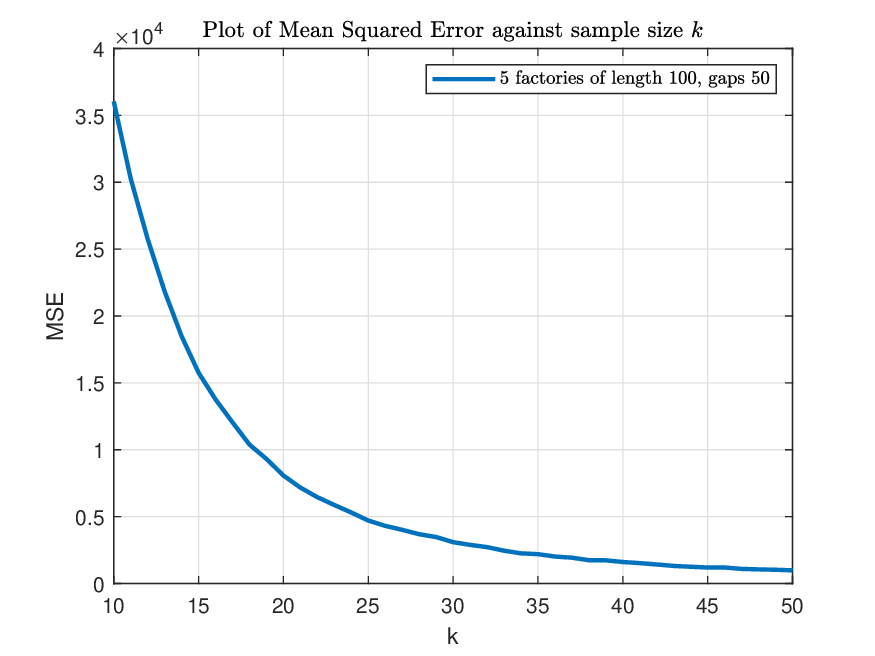}}\\
\scalebox{.45}{\includegraphics{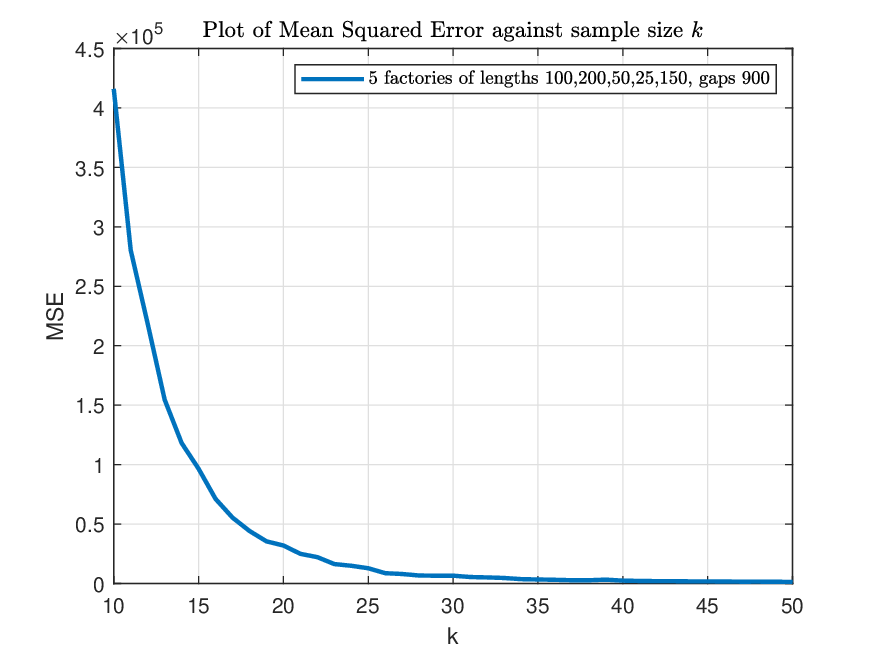}}\ \scalebox{.45}{\includegraphics{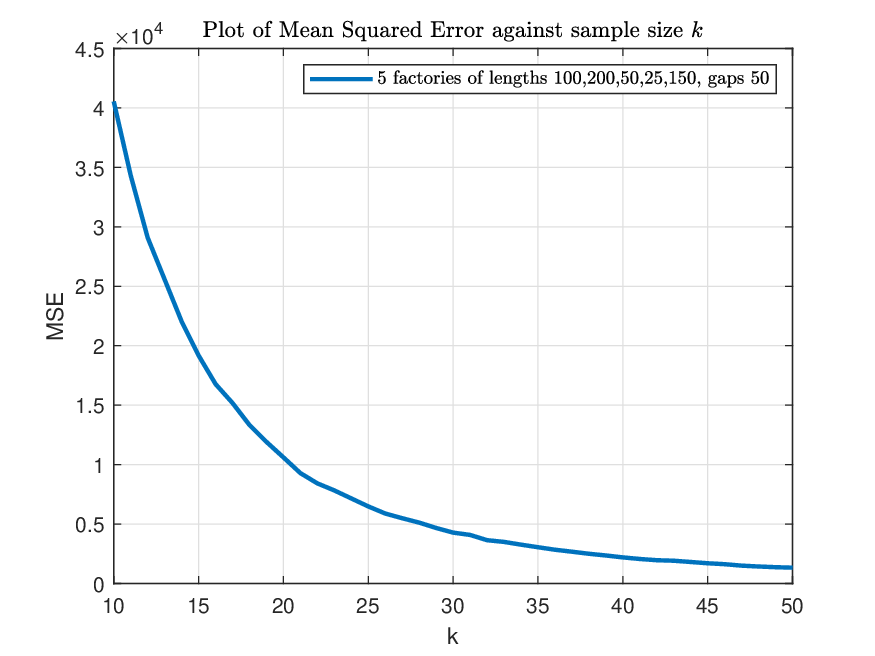}}\\
\scalebox{.45}{\includegraphics{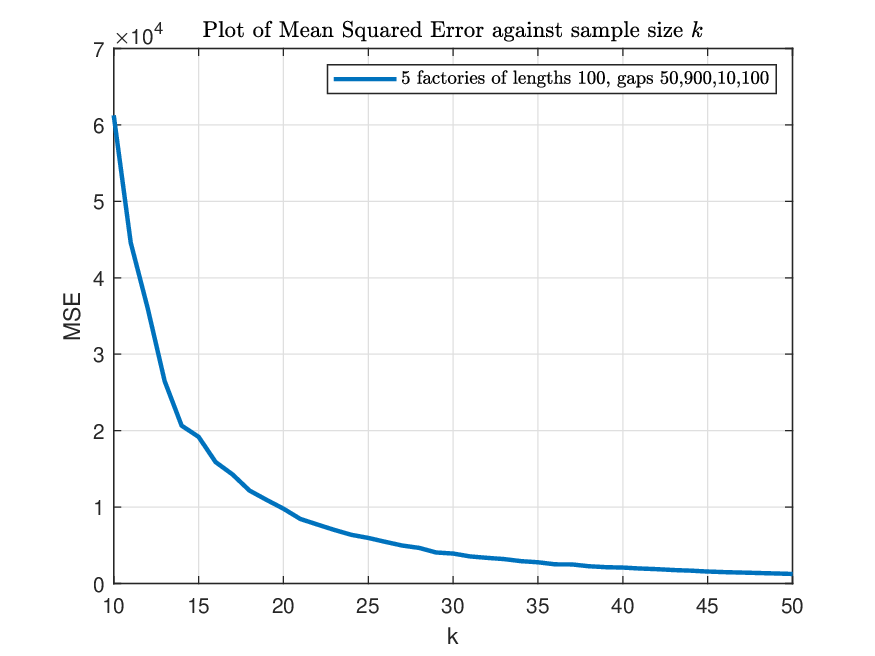}}\ \scalebox{.45}{\includegraphics{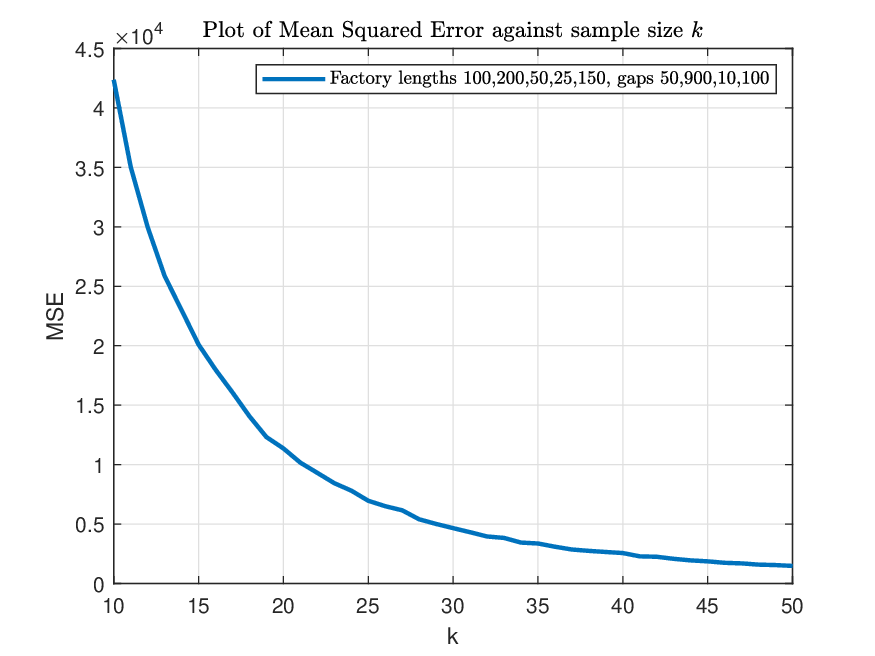}}
\caption{\label{figure:MSE} Mean Squared Error of 10,000 estimations by the MFP plotted against sample size $k$.}
\end{center}
\end{figure}

As detailed in previous sections, we have designed our MFP estimator to be robust when $k$ is sufficiently large compared to $l$. In Figure \ref{figure:MSE}, we can see that when $l=5$, the estimator is terrible when $k=10$, but improves dramatically with just a few additional samples, before showing more modest improvement after $k=30$.

Figure \ref{figure:MSE} shows the estimator performing well in a range of situations. When gaps are large relative to the factory sizes (from top, left to right, see the first and third plots) the estimator performs well as the chance of the four largest gaps between the ordered samples splitting the samples into their respective factories is very high. If this is the casee, all error comes down to the variance within the GTP and GTP-UM, which we proved to be the least possible for an unbiased estimator in Section \ref{sec:mvue}.

However, when gaps are small relative to factory size (see the second and fourth plots) we also get good performance, once our sample size is sufficiently large. What might seem like a harder problem is ultimately the same - when there are enough samples, two ordered samples from the same factory are still likely to be closer together than two samples from different factories.

Interestingly, we get noticeably \emph{better} performance for small $k$. The major reason this happens is that when we are missing the first factory entirely, smaller gaps mean our smallest sample is still going to be relatively small. This mitigates the extent to which we end up overestimating the size of the first factory.

In fact, this is so significant that the Mean Squared Error in the first plot is larger than $500^2$ for $k=10$. This is an enormous error as the total factory size is 500. This happens because we often estimate the first factory to have a size of over 1000 if there was no sample from the true first factory.

Another reason smaller gaps are easier is when the gaps are small, mistaking a gap between two samples from the same factory as those from different factories is less critical an error, since the true gap size is unlikely to be much bigger.

For completeness, in the fifth and sixth plots we demonstrate similarly good performance when taking gaps of different lengths. There is no reason to have thought our estimator would do worse in this situation. But if we knew that the gaps were the same length, this reduces the number of gap-related unknowns from $l$ to $1$. This is a considerably different problem that we suspect can be approached in a way that is effective even if we do not have a sample from every factory, making it more robust than our method when the number of samples are low.

In the following section, we explore what progress can be made when we assume our gaps are both \emph{known} and \emph{fixed} using an adapted version of the approach seen in the original GTP.

\section{Restricting to Equal Factory Size and Fixed, Known Gaps}
\label{sec:restriction}

In this section we are estimating the equal size of each factory, $N$, given $k$ samples with $G$ the fixed, known gap size, and $l$ the number of factories. Let $x_{\parentheses{1}},\ldots,x_{\parentheses{k}}$ be the observed ordered serial numbers, with the maximum $x_{\parentheses{k}} \ = \ M$.

For a small number of samples, we can find an unbiased estimator by looking at $\expectation{M}$.  Directly looking at the probability mass function given $k$ samples gives

\begin{equation}
    \probability{M = x} \ = \ \frac{\binom{x  - G\floor{\frac{x}{N+G}} -1}{k - 1}}{\binom{l N}{k}}
\end{equation}

\noindent for $x \in T$. This however is not easy to use or manipulate. Instead, we can look at the value of the serial numbers differently, as a scaled up version of the original one factory problem, since we are only looking at the expected value, and not the distribution. We shall split $M = \Tilde{x} + H$. Here $\Tilde{x} = M - G\floor{\frac{M}{N+G}}$, i.e., it corresponds to what the serial number $M$ would be if there were no gaps. Therefore,

\begin{equation}
    \probability{\Tilde{x} = x} \ = \ \frac{\binom{x-1}{k-1}}{\binom{l N}{k}}
\end{equation}

\noindent which has expectation (see \cite{LM})
\begin{equation}
    \expectation{\Tilde{x}} \ = \ \frac{k}{k+1} \parentheses{l N + 1}.
\end{equation}

Looking at $H$, where $H \ = \ G\floor{\frac{M}{N+G}}$, yields

\begin{equation}
    \expectation{H} \ = \ \sum_{x = k}^{lN} \frac{\binom{x-1}{k-1}}{\binom{lN}{k}} \cdot G\floor{\frac{x}{N}}.
\end{equation}

Let us take the case of $k \leq N$. We can re-index the sum and use the hockey stick identity \eqref{eq:hockeystickidentity} for sums of binomial coefficients to get

\begin{equation}\label{eq:H formula}
\begin{split}
    \expectation{H} \ &= \ \frac{G}{\binom{lN}{k}}\sum_{t = 1}^{l-1} \sum_{x = tN + 1}^{lN} \binom{x-1}{k-1} \\
    &= \ \frac{G}{\binom{lN}{k}}\sum_{t = 1}^{l-1} \parentheses{ \sum_{x = k}^{lN} \binom{x-1}{k-1} -  \sum_{x = k}^{tN} \binom{x-1}{k-1}} \\
    &= \ \frac{G}{\binom{lN}{k}} \sum_{t = 1}^{l-1} \parentheses{\binom{lN}{k} - \binom{tN}{k}} \\
    &= \ G\parentheses{l-1} - \frac{G}{\binom{lN}{k}}\sum_{t = 1}^{l-1} \binom{tN}{k} \\
    &= \ Gl - \frac{G}{\binom{lN}{k}}\sum_{t = 1}^{l} \binom{tN}{k}.
\end{split}
\end{equation}

This is very useful as for small values of $k$ a closed form can be directly calculated, and for larger values, we can take an approximation as $N \to \infty$. Taking $k=1$ gives 

\begin{equation}
\begin{split}
    \frac{G}{\binom{lN}{1}}\sum_{t = 1}^{l} \binom{tN}{1} \ &= \ \frac{G}{lN} \sum_{t=1}^{l} tN \\
    &= \ \frac{GN}{lN} \sum_{t=1}^{l} t \\
    &= \ \frac{GlN\parentheses{l+1}}{2lN} \\
    &= \ \frac{G\parentheses{l+1}}{2}.
\end{split}
\end{equation}

Substituting this back yields,
\begin{equation}
\begin{split}
    \expectation{M} \ &= \ \expectation{\Tilde{x}} + \expectation{H} \\
    &= \ \frac{lN + 1}{2} + Gl - \frac{G\parentheses{l+1}}{2} \\
    &= \ \frac{l}{2} N + \frac{G\parentheses{l-1} + 1}{2}
\end{split}
\end{equation}

\noindent and rearranging this, we obtain

\begin{equation}
    \hat{N} \ = \ \frac{1}{l} \parentheses{2M - G\parentheses{l-1} - 1}.
\end{equation}

\begin{rek}
    We remark that we can find such a closed form for $k < 5$ as there is a formula for the roots of a $k^\text{th}$ degree polynomial in this case. 
\end{rek}

For general $k$, we can approximate $\hat{N}$. The only issue arises in calculating the final sum in \eqref{eq:H formula}. For this, we can expand the binomial coefficient and then use the two leading terms from Faulhaber's formula for the sum of $k^{\text{th}}$ powers of natural numbers (see \cite{CG} for derivation and more details) to approximate the sum for large $N$. For the following, take $k \ll lN$:

\begin{equation}
\begin{split}
    \frac{G}{\binom{lN}{k}}\sum_{t = 1}^{l} \binom{tN}{k} \ &= \ \frac{G}{k! \binom{lN}{k}}\sum_{t = 1}^{l} tN \parentheses{tN-1} \cdots \parentheses{tN-\parentheses{k-1}} \\
    &\approx \ \frac{G}{\parentheses{lN}^{k}}\sum_{t = 1}^{l} \parentheses{ t^k N^k } \\
    &\approx \ \frac{GN^k}{\parentheses{lN}^{k}} \parentheses{\frac{1}{k+1} l^{k+1} + \frac{1}{2} l^k} \\
    &= \ G \parentheses{\frac{1}{k+1}l + \frac{1}{2}}. \\
\end{split}
\end{equation}

Substituting this back into \eqref{eq:H formula}

\begin{equation} \label{eq: Hexpectedvalue}
    \expectation{H} \ \approx \ G\parentheses{\frac{k}{k+1} l - \frac{1}{2}}
\end{equation}

\noindent and so, we get 

\begin{equation} \label{eq: Mexpectedvalue}
     \expectation{M} \ \approx \ \frac{k}{k+1} \parentheses{lN+1} + G\parentheses{\frac{k}{k+1}l - \frac{1}{2}}.
\end{equation}

From this we can rearrange to $N$ to find an approximate estimator $\hat{N}$:
\begin{equation}
\begin{split}
    \hat{N} \ &\approx \ \frac{1}{l} \parentheses{\frac{k+1}{k} M - \frac{k+1}{k}G\parentheses{\frac{k}{k+1}l - \frac{1}{2}} - 1} \\
    &\approx \ \frac{1}{l} \parentheses{\frac{k+1}{k} M - Gl + G\frac{k+1}{2k} - 1}.
\end{split}
\end{equation}


As detailed previously, we expect that the approximate estimator should work well for $k \ll lN$. Further, we expect nicer behaviour for small values of $G$ since, in \eqref{eq: Mexpectedvalue}, the first term is independent of $G$, whereas the second term which corresponds to $\expectation{H}$ depends on it. Therefore, any error in our approximation is exaggerated for $G$ large.   

\begin{figure}[h]
\begin{center}
\scalebox{.44}{\includegraphics{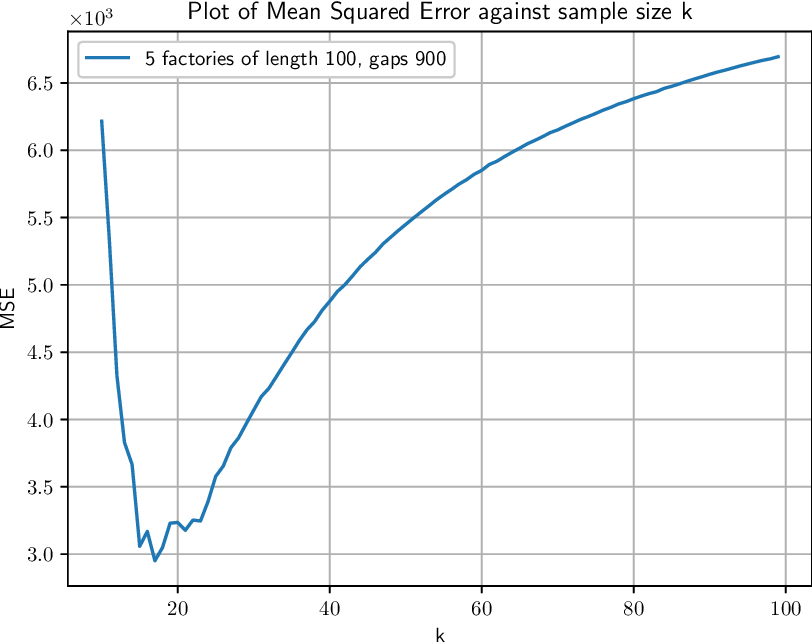}}\ \scalebox{.44}{\includegraphics{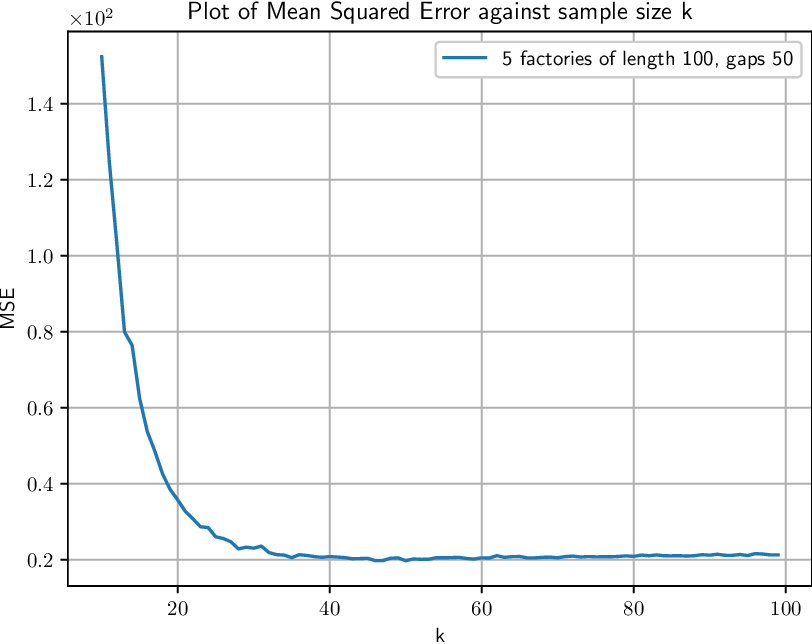}} \\
\scalebox{.44}{\includegraphics{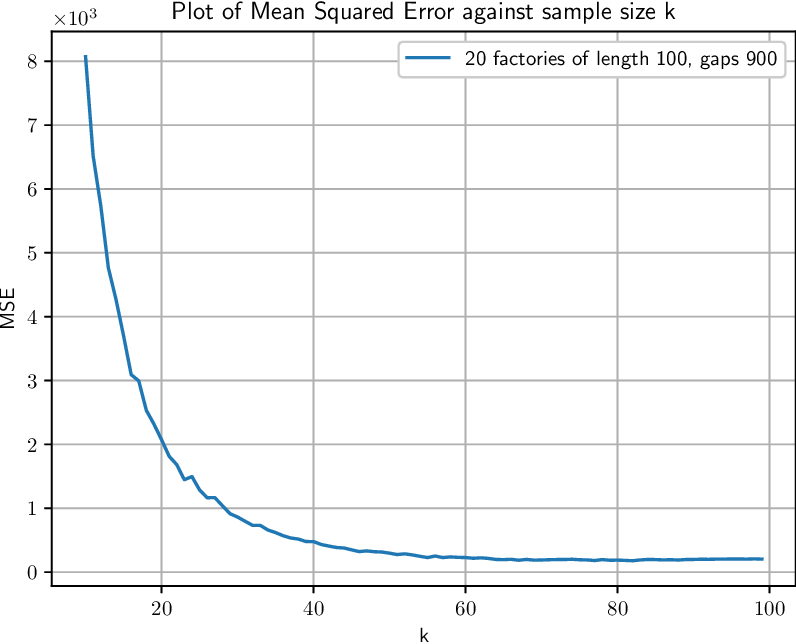}}\ \scalebox{.44}{\includegraphics{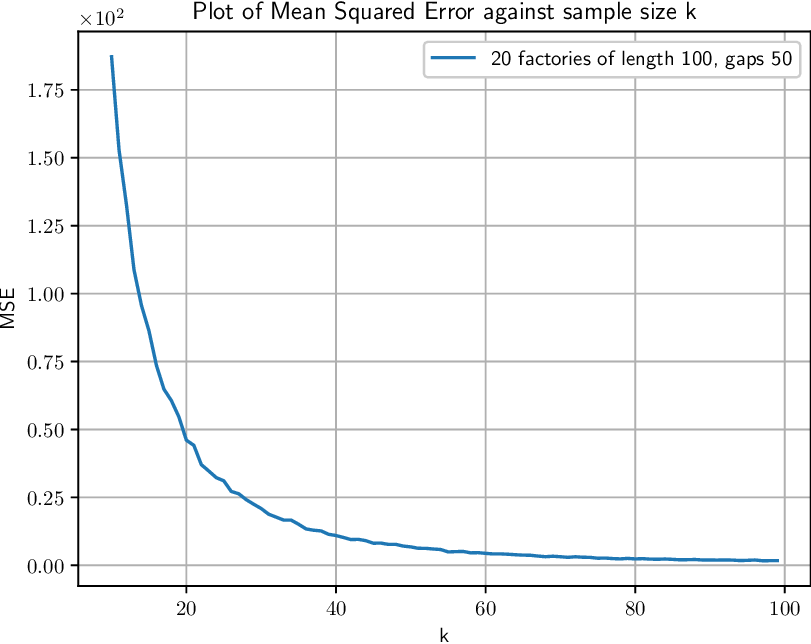}}
\caption{\label{figure:MSE2} Mean Squared Error of 10,000 estimations by the MFP plotted against sample size $k$.}
\end{center}
\end{figure}

Figure \ref{figure:MSE2} shows that the estimator has lower MSE than that for varied gaps and factory lengths (see scale on the graphs). As expected, when the number of factories is larger, our MSE is eventually smaller, as can be seen in the graphs showing the cases for $l = 20$, compared to those for $l = 5$.

Further, in the top left graph, showing 5 factories of length 100 with gaps of size 900, we can see the balancing of the two main effects from changing $k$: first, as $k$ increases, we are more likely to get a sample from every factory and the observed maximum $M$ is on average closer to $\expectation{M}$. On the other hand, the approximation in \eqref{eq: Hexpectedvalue} becomes less precise, and this eventually overpowers the former. For less than 25 samples, the MSE is still orders of magnitude smaller than that of the MFP estimator with unknown gap size and factory length.  

\begin{figure}[h]
\begin{center}
\scalebox{.44}{\includegraphics{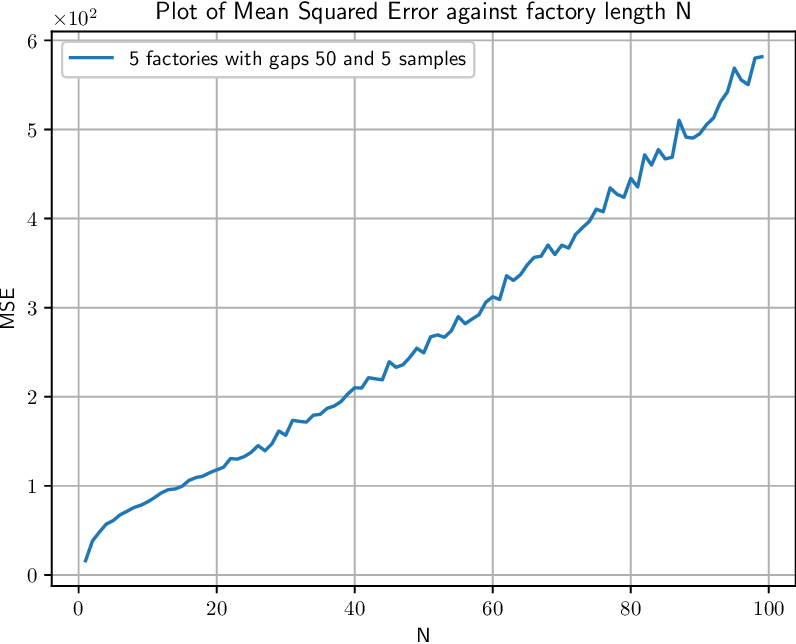}}\ \scalebox{.44}{\includegraphics{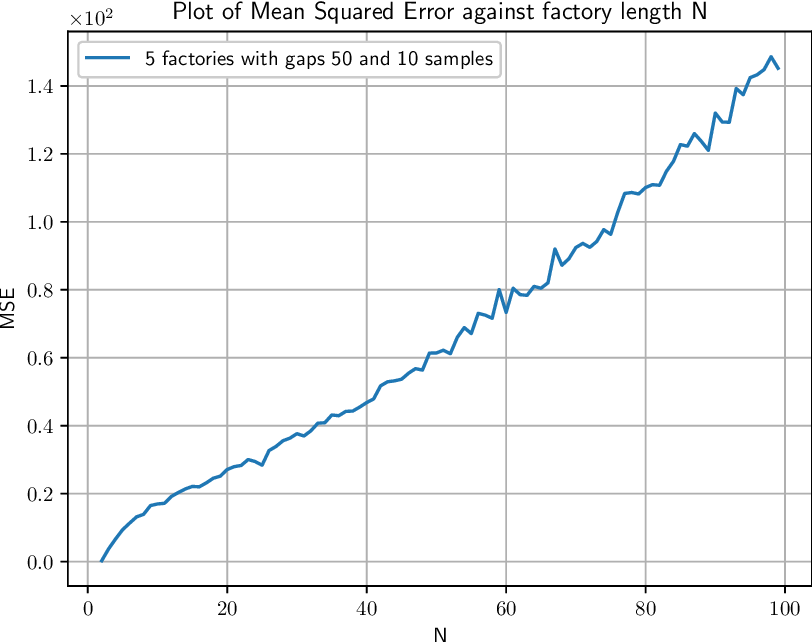}}\\
\scalebox{.44}{\includegraphics{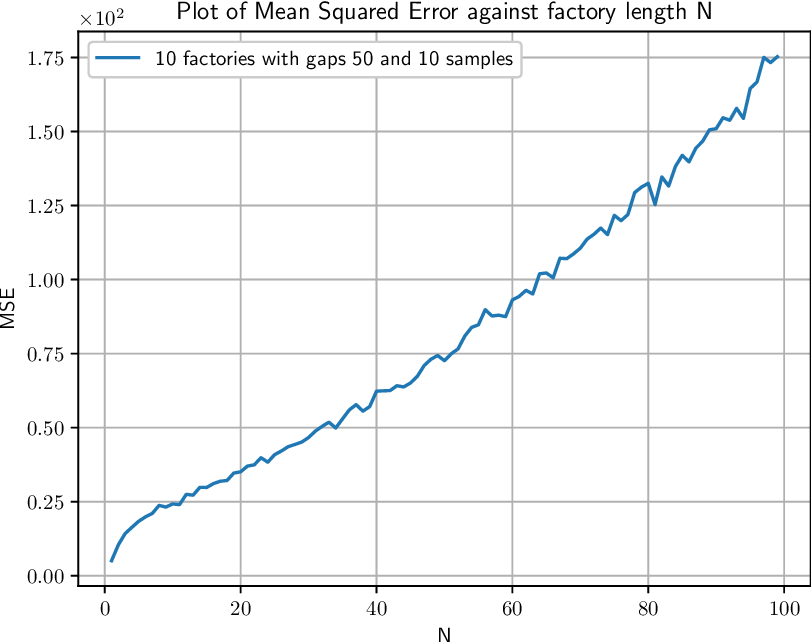}}\ \scalebox{.44}{\includegraphics{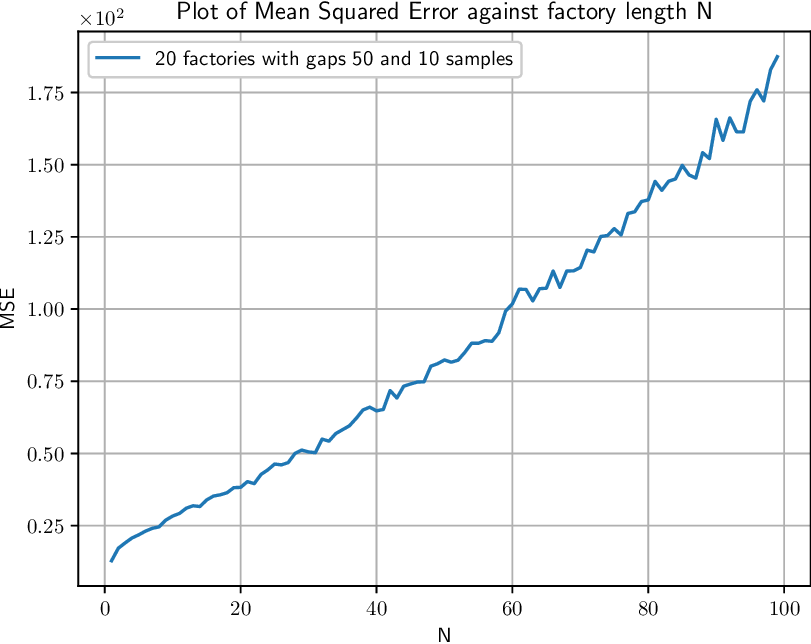}}
\caption{\label{figure:MSE_changeN_fixk} Mean Squared Error of 10,000 estimations by the MFP plotted against factory length $N$ for varied choices of fixed sample size $k$ and factory number $l$.}
\end{center}
\end{figure}

\begin{figure}[h]
\begin{center}
\scalebox{.44}{\includegraphics{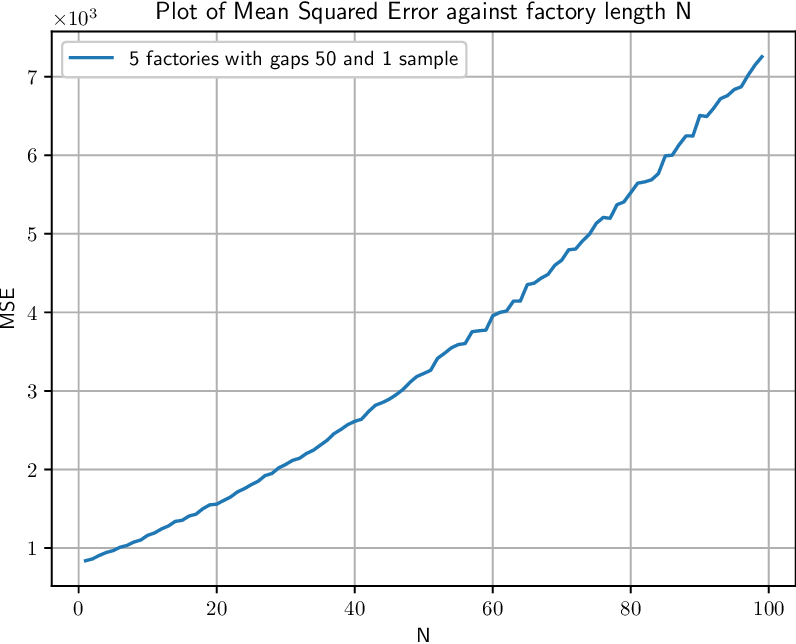}}\
\scalebox{.44}{\includegraphics{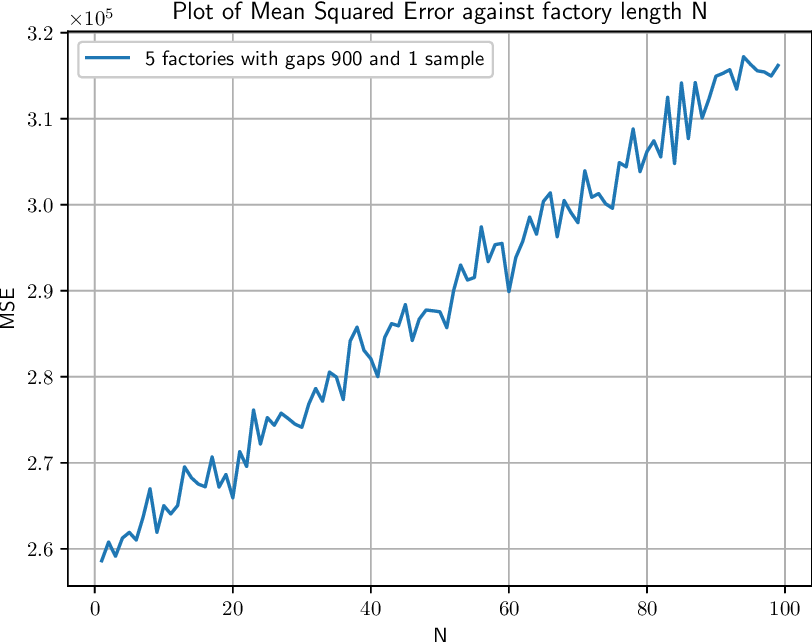}}
\caption{\label{figure:MSE_changeN_fixk=1} Mean Squared Error of 10,000 estimations by the MFP plotted against factory length $N$ for fixed sample size $k = 1$.}
\end{center}
\end{figure}

Varying the factory length, as can be seen in Figure \ref{figure:MSE_changeN_fixk}, we can see that in all shown cases, the MSE does not rise above $10^3$ for $N \leq 100$. Further, as the number of factories increases, our MSE does not change by much.

Although there is the unbiased exact estimator in the case $k = 1$ in Figure \ref{figure:MSE_changeN_fixk=1}, the MSE is expectionally high. As expected, it is not possible to reliably estimate the factory length well from a single sample, since we have such a high variance when we draw the sample.

\section{Extensions}
\label{sec:ext}

Considering that the Multiple Factories Problem is a novel problem, there is plenty to explore and we encourage future research to build on the foundations laid in this paper. The following is a list of open problems that we considered whilst writing:

\begin{itemize}
    \item Can we make any progress if we make the number of factories $l$ unknown?
    \item What happens to the probability of missing a factory when we take $k,l \to \infty$ and keep $N$ finite?
    \item What significant improvements can be made to our method to attack the full MFP? Could we possibly find an unbiased estimator for the total number of tanks?
    \item What is the variance of our approach to the MFP when we have fixed lengths, and fixed and known gaps?
    \item How much progress can we make when we know the gaps are fixed, but not known? What performance do we get from estimating the fixed gap and plugging it in?
\end{itemize}

\section{Conclusion}
\label{sec:conc}

Returning to the real world, we can conclude that setting large gaps between the serial numbers of different factories is beneficial in hiding your total production. Consider the following real situation: you are given a sample of captured tanks, where the smallest serial number is already a large number. Perhaps the first factory is very productive, or perhaps all of our captured tanks are just from later factories. If we answer this incorrectly, the error in our overall estimation of productivity will increase when we have larger gaps. However, our approach remains robust in the situation with large gaps, due to the very high probability that we identify the gaps between factories correctly given enough samples. 

In the situation with small gaps, our approach will always do well. It is significant that our approach to estimating total production is robust against any gaps between factories you can choose - large or small, fixed or widely varied. Ultimately, it does not matter whether we correctly identify which tank is from which factory, as long as our estimate of the total number of tanks is close.

When the factory length is fixed and the gap size is fixed and known, we can make considerably more progress in making a lower variance estimate of the total number of tanks. This is true even when the sample size is smaller than the number of factories, at which point we make several orders of magnitude improvement in mean squared error.

We conclude that the Multiple Factories Problem is considerably harder than the German Tank Problem, but still tractable provided we have enough samples, and very tractable if we are given information about the gaps between the factories.

\appendix


\section{Proofs of Identities}
\label{sec:identity}

\noindent Identity I: For all $N \geq k$:
\begin{equation}
    \sum_{m=k-b+1}^{N-b+1} m\frac{\binom{m-1}{k-b} \binom{N-m}{b-1}}{\binom{N}{k}} \ = \ \frac{\parentheses{N+1}\parentheses{k-b+1}}{\parentheses{k+1}}
\end{equation}

and \\

\noindent Identity II: For all $N \geq k$:
\begin{equation}
\begin{split}
    \sum_{m=k-b+1}^{N-b+1} m^2\frac{\binom{m-1}{k-b} \binom{N-m}{b-1}}{\binom{N}{k}} \ & = \ \frac{(k-b+1)(k-b+2)(N+2)(N+1)}{(k+2)(k+1)} \\ &- \  \frac{\parentheses{N+1}\parentheses{k-b+1}}{\parentheses{k+1}}.
\end{split}
\end{equation}

The proofs require another identity, Identity III:
\begin{equation} \label{eq:identityIII}
    \binom{a+b+k+1}{a+b+1} \ = \ \sum_{i=0}^k \binom{a+i}{a} \binom{b+k-i}{b}.   
\end{equation}

\begin{proof} (of Identity III)

We use proof by induction on $b$ to prove the identity and thus by symmetry of $a$ and $b$, it holds for all $a,b$. Consider the base case where $a=b=0$:

\begin{align}
    \binom{k+1}{1} \ &= \ k + 1 \\
    &= \ \sum_{i=0}^k \binom{i}{0}\binom{k-i}{0}.
\end{align}

This is true for all $k$. Now, suppose that the identity holds for $a,b$, then we shall sum over the left hand side of the identity from $k=0$ to $m \in \mathbb{N}$:

\begin{align}
    \sum_{k=0}^m \binom{a+b+k+1}{a+b+1} \ &= \ \binom{a+b+m+2}{a+b+2} \\
    &= \ \binom{a + \parentheses{b+1} + m + 1}{a + \parentheses{b+1} + 1},
\end{align}

where we used the hockey stick identity:

\begin{equation}\label{eq:hockeystickidentity}
    \sum_{i=r}^n \binom{i}{r} \ = \ \binom{n+1}{r+1}.
\end{equation}

Summing over the right hand side, re-indexing the summation and then applying the hockey stick identity gives

\begin{align*}
    \sum_{k=0}^m \sum_{i=0}^k \binom{a+i}{a} \binom{b+k-i}{b} \ &= \ \sum_{i=0}^m \sum_{k=i}^m  \binom{a+i}{a} \binom{b+k-i}{b} \\
    &= \ \sum_{i=0}^m \binom{a+i}{a} \sum_{k=i}^m \binom{b+k-i}{b} \\
    &= \ \sum_{i=0}^m \binom{a+i}{a} \binom{\parentheses{b+1}+k-i}{b+1} \\
    &= \ \binom{a + \parentheses{b+1} + m + 1}{a + \parentheses{b+1} + 1}. \numberthis
\end{align*}\end{proof}

Now, we can quickly prove Identities I and II by applying Identity III.

\begin{proof} (of Identity I). We have
\begin{align*}
     &\sum_{m=k-a+1}^{N-a+1} m\frac{\binom{m-1}{k-a} \binom{N-m}{a-1}}{\binom{N}{k}}\\
     = \ &\binom{N}{k}^{-1} \sum_{m=k-a+1}^{N-a+1} m \cdot \binom{m-1}{k-a} \binom{N-m}{a-1}\\
     = \ &\binom{N}{k}^{-1} \sum_{m=k-a+1}^{N-a+1} \frac{m!}{\parentheses{k-a}!\parentheses{m-k+a-1}!} \cdot \frac{\parentheses{N-m}!}{\parentheses{a-1}!\parentheses{N-m-a+1}!} \\
     = \ &\frac{k-a+1}{\binom{N}{k}} \cdot \sum_{m=k-a+1}^{N-a+1} \binom{m}{k-a+1} \cdot \binom{N-m}{a-1}\\
     = \ &\frac{\parentheses{N+1}\parentheses{k-a+1}}{\parentheses{k+1}}. \numberthis
\end{align*}
\end{proof}

\begin{proof} (of Identity II). We have
\begin{align*}
    \sum_{m=k-a+1}^{N-a+1} m^2 \frac{\binom{m-1}{k-a}\binom{N-m}{a-1}}{\binom{N}{k}} \ &= \ \binom{N}{k}^{-1}\sum_{m=k-a+1}^{N-a+1} \parentheses{m+1}m\binom{m-1}{k-a}\binom{N-m}{a-1} \\
    & \hspace{0.7cm} - \ \binom{N}{k}^{-1}\sum_{m=k-a+1}^{N-a+1} m\binom{m-1}{k-a}\binom{N-m}{a-1}\\
    &= \ \binom{N}{k}^{-1}\sum_{m=k-a+1}^{N-a+1} \frac{\parentheses{m+1}!}{\parentheses{k-a}!\parentheses{m-k+a-1}!}\binom{N-m}{a-1} \\
    & \hspace{0.7cm} - \ \frac{\parentheses{N+1}\parentheses{k-a+1}}{\parentheses{k+1}} \\\
    &= \ \frac{\parentheses{k-a+1}\parentheses{k-a+2}\parentheses{N+2}\parentheses{N+1}}{\parentheses{k+2}\parentheses{k+1}} \\
    & \hspace{0.7cm} - \ \frac{\parentheses{N+1}\parentheses{k-a+1}}{\parentheses{k+1}}. \numberthis
\end{align*}
\end{proof}

\section{Code}
\subsection{First Approach at MFP}
\label{sec:codeandy}

These are MATLAB programs and require the standard Statistics and Machine Learning Toolbox.

\begin{itemize}
\item A function that generates a sample uniformly without replacement from some factories.

\begin{lstlisting}[ basicstyle = \ttfamily\tiny]
function sample = genSample(factories,r)
% generate r uniform w/o replacement from factories =
% [a_1 a_2 a_3 ... a_2m]
m = length(factories)/2;
N = 0;
for i = 1:m % get size N first
    N = N + factories(2*i) - factories(2*i-1) + 1;
end
source = zeros(1,N);
pos = 1;
for i = 1:m % build source to sample from
    fact = factories(2*i-1):1:factories(2*i);
    source(pos:pos+length(fact)-1) = fact;
    pos = pos + length(fact);
end
sample = sort(randsample(source,r));
end
\end{lstlisting}

\item A function that performs the MFP estimate on sample

\begin{lstlisting}[ basicstyle = \ttfamily\tiny]
function estN = advMultGerman(sample,l,lowerKnown)
% our new german tank with l factories
% output is of form a_1,a_2,..,a_2k where factories are from these ranges
% handles when we have a singleton sample of a factory, 
if ~exist('lowerKnown','var')
    % third parameter does not exist, so default it to something
    lowerKnown = true;
end
r = length(sample);
% sort gaps
gaps = zeros(1,r-1);
for i = 1:r-1
    gaps(i) = sample(i+1) - sample(i);
end
[~,bigGaps] = maxk(gaps,l-1); % find indices where l-1 largest gaps appear
bigGaps = sort(bigGaps);
% disp(bigGaps)
% compute estimate
estN = bruteAdv(sample,l,lowerKnown,bigGaps,r);
end
\end{lstlisting}

\item An auxiliary function for the above.

\begin{lstlisting}[ basicstyle = \ttfamily\tiny]
function estN = bruteAdv(sample,l,lowerKnown,bigGaps,r)
% we handle singletons bby
Ns = zeros(1,l); % number in each factory stored here
bad = []; % store bad subsample indexes here
Nsum = 0;
subsample = sample(1:bigGaps(1));
ksum = length(subsample);
if lowerKnown == true
    % one side on first gap
    Ns(1) = german(subsample);
    Nsum = Nsum + Ns(1);
else % two sided on first gap
    Ns(1) = twoSideGerman(subsample);
    Nsum = Nsum + Ns(1);
end
for i = 2:l-1
    subsample = sample(bigGaps(i-1) + 1:bigGaps(i));
    k = length(subsample);
    if k == 1
        % add this bad subsample to index set
        bad = [bad,i];
    else
        ksum = ksum + k;
        Ns(i) = twoSideGerman(subsample);
        Nsum = Nsum + Ns(i);
    end
end
subsample = sample(bigGaps(l-1)+1:r);
k = length(subsample);
if k == 1
    bad = [bad,l];
else
    ksum = ksum + k;
    Ns(l) = twoSideGerman(subsample);
    Nsum = Nsum + Ns(l);
end
Nhat = Nsum/ksum;
% handle bad subsamples
for element = bad
    Ns(element) = Nhat;
end
% disp(Ns)
estN = sum(Ns);
end
\end{lstlisting}

\item An auxiliary function that performs the original GTP on a sample.

\begin{lstlisting}[ basicstyle = \ttfamily\tiny]
function estN = german(sample)
% classic German Tank estimate for N, sample size r and max m
r = length(sample);
m = max(sample);
estN = m*(1+1/r)-1;
end
\end{lstlisting}

\item An auxiliary function that performs the two-tailed GTP-UM on a sample.

\begin{lstlisting}[ basicstyle = \ttfamily\tiny]
function estN = twoSideGerman(sample)
% german tank estimator where lower end is unknown
% estimate for N = N2-N1+1, [N1, N2], sample size r and spread s
% DOES NOT WORK WHEN HAVE ONLY ONE SAMPLE
r = length(sample);
s = max(sample) - min(sample);
estN = s*(1+2/(r-1))-1;
end
\end{lstlisting}

\item A function that finds the MSE of the MFP.

\begin{lstlisting}[ basicstyle = \ttfamily\tiny]
function [avg, mse] = simulMFP(factories, trials, k, l, N)
% simulate sample mean and MSE of MFP for given factories, samples, #
% factories, and # repeats trials, and total tanks N
results = zeros(1,trials);
ideal = N + results;
for i = 1:trials
    sample = genSample(factories,k);
    results(i) = advMultGerman(sample,l);
end
avg = mean(results);
squaredError = (ideal - results).^2;
mse = mean(squaredError);
end
\end{lstlisting}

\end{itemize}

\subsection{Fixed Lengths and Fixed, Known Gaps}
\label{sec:codekishan}

These are Python programs used in Jupyter notebook.

\begin{itemize}
\item Importing modules

\begin{lstlisting}[ basicstyle = \ttfamily\tiny]
import random
import matplotlib.pyplot as plt
import numpy as np
import pandas as pd
plt.rcParams['text.usetex'] = True
\end{lstlisting}

\item Function to generate samples
\begin{lstlisting}[ basicstyle = \ttfamily\tiny]
def samples(k,g,l,N):
    list_of_tanks = []
    for num in range(l):
        list_of_tanks += list(range(num*(N+g)+1,num*(N+g) + N + 1))
    list_of_samples = random.sample(list_of_tanks,k=k)
    list_of_samples.sort()
    return list_of_samples
\end{lstlisting}

\item Exact estimator for $k = 1$ case
\begin{lstlisting}[ basicstyle = \ttfamily\tiny]
def mle1(g,l,N):
    max_tank = samples(k,g,l,N)[-1]
    return 1/l * ( 2*max_tank - (l-1)*g - 1)

def mle1_repeat(number_of_tests,gap,number_of_factories,length):
    data = []
    for i in range(number_of_tests):
        data.append(mle1(gap,number_of_factories,length))
    return[np.mean(data),np.mean([(datum - length)**2 for datum in data])]

\end{lstlisting}

\item Approximate estimator for general $k$
\begin{lstlisting}[ basicstyle = \ttfamily\tiny]
def mle_estimator(k,g,l,N):
    #given fixed k,g,l, we can estimate N 
    max_tank = samples(k,g,l,N)[-1]
    N_hat = 1/l * ((k+1)*max_tank/k - g*l + g*(k+1)/(2*k) - 1)
    return N_hat

def mle_estimator_repeat(number_of_tests,samples,gap,number_of_factories,length):
    data = []
    for i in range(number_of_tests):
        data.append(mle_estimator(samples,gap,number_of_factories,length))
    return[np.mean(data),np.mean([(datum - length)**2 for datum in data])]
\end{lstlisting}

\item Used to generate graphs changing $k$
\begin{lstlisting}[ basicstyle = \ttfamily\tiny]
g = 50
l = 20
N = 100
A = 1
c = 0.7
minimum = 10
maximum = 50
tests = 10000
data =[[],[]]
our_range = range(minimum,maximum)


for k in our_range:
    current_data = mle_estimator_repeat(tests,k,g,l,N)
    data[0].append(current_data[0])
    data[1].append(current_data[1])

plt.plot(our_range,data[1],'tab:blue')
plt.legend([str(l)+' factories of length '+str(N)+', gaps '+str(g)])

plt.xlabel('k')
plt.ylabel('MSE')
plt.title('Plot of Mean Squared Error against sample size \\textit{\emph{k}}')
plt.ticklabel_format(style='sci',axis='y',scilimits=(0,0))

plt.grid()

plt.savefig(str(l)+'l'+str(N)+'N'+str(g)+'g_graph.pdf')

\end{lstlisting}

\item Used to generate graphs for fixed $k$
\begin{lstlisting}[ basicstyle = \ttfamily\tiny]
g = 50
l = 20
k = 10
A = 1
c = 0.7
maximum = 100
tests = 10000
data =[[],[]]
our_range = range(int(np.ceil(k/l)),maximum)


for N in our_range:
    current_data = mle_estimator_repeat(tests,k,g,l,N)
    data[0].append(current_data[0])
    data[1].append(current_data[1])

plt.plot(our_range,data[1],'tab:blue')
plt.legend([str(l)+' factories with gaps '+str(g)+' and ' + str(k) + ' samples'])

plt.xlabel('N')
plt.ylabel('MSE')
plt.title('Plot of Mean Squared Error against factory length \\textit{\emph{N}}')
plt.ticklabel_format(style='sci',axis='y',scilimits=(0,0))

plt.grid()

plt.savefig(str(l)+'l'+str(k)+'k'+str(g)+'g_graph.pdf')
\end{lstlisting}

\end{itemize}

\ \\

\end{document}